# FEASIBLE LOW-THRUST TRAJECTORY IDENTIFICATION VIA A DEEP NEURAL NETWORK CLASSIFIER

**Ruida Xie[*], Andrew G. Dempster [†]**

In recent years, deep learning techniques have been introduced into the field of trajectory optimization to improve convergence and speed. Training such models requires large trajectory datasets. However, the convergence of low thrust (LT) optimizations is unpredictable before the optimization process ends. For randomly initialized low thrust transfer data generation, most of the computation power will be wasted on optimizing infeasible low thrust transfers, which leads to an inefficient data generation process. This work proposes a deep neural network (DNN) classifier to accurately identify feasible LT transfer prior to the optimization process. The DNN-classifier achieves an overall accuracy of 97.9%, which has the best performance among the tested algorithms. The accurate low-thrust trajectory feasibility identification can avoid optimization on undesired samples, so that the majority of the optimized samples are LT trajectories that converge. This technique enables efficient dataset generation for different mission scenarios with different spacecraft configurations.

**INTRODUCTION**

Low thrust trajectory optimization can be formulated as a constrained Optimal Control Problem (OCP) and can be solved by indirect and direct methods[1]. The drawbacks of these methods include relying heavily on good initial guesses[2-4], no guaranteed convergence[5, 6], and the long computational time[7]. On the other hand, in recent years, machine learning methods are being studied and introduced into this research field to replace the conventional methods[1]. The roadmap is to train regression models as surrogates to approximate the transfer cost, such that the complex optimization process can be avoided. This approach has been applied to solve complex multi-target visiting sequence design in the Global Trajectory Optimization Competitions (GTOCs) [8-10], landing guidance and control[6, 11], etc.

During the training of such models, one of the major obstacles is to generate a high-quality low-thrust (LT) trajectory database efficiently for different mission scenarios with different spacecraft configurations. The current methods, which repeatedly perform optimization for random selected target pairs and transfer parameters, suffers from low convergence rate. The feasibility of LT transfers is unpredictable before the optimization process ends. In this work, by the feasibility of a LT transfer, we mean the convergence of that LT trajectory optimization. In our testing case, the convergence rate is approximately 7% if the targets are selected completely randomly.

---

[*] Ph.D. candidate, Australian Centre for Space Engineering Research (ACSER), University of New South Wales (UNSW), Sydney, NSW, 2052, Australia

[†] Professor, Australian Centre for Space Engineering Research (ACSER), University of New South Wales (UNSW), Sydney, NSW, 2052, Australia

The current solutions include limiting the difference between some orbital parameters of asteroids pairs[2, 3] or constructing phase indicators to rank the low-thrust (LT) transfer convergence possibility[12, 13]. However, these methods have limited performance as they can only be used as a reference. Furthermore, using these methods may filter out some feasible transfers, which leads to a biased training database.

Motivated by the above observations, this work provides the methodology for training classification models for accurate low-thrust (LT) transfer feasibility identification, so that the convergence of the low-thrust optimization can be predicted. The training data is generated by using Sims-Flanagan (SF) model[14], and features are constructed according to the parameters required in the original LT optimization problems. An embedded feature importance estimator is presented for preliminary feature selection. Eventually, the DNN-classifier is constructed and the hyperparameters is optimized by parallel Bayesian Optimization process.

In addition, we compared the DNN-classifier performance with other four representative classification models, including the Ensemble learning model (Ensembled trees), K-Nearest Neighbors (KNN) model, Support Vector Machine (SVM) and Tree based model. The training set size requirements, data augmentation and optimization process of classifiers' hyperparameters are presented in the simulation section. The DNN-classifier has the most advanced performances on all evaluation metrics.

Although the purpose of the proposed DNN-classifier is to increase the convergence during LT data generation, the proposed method can also be used as an accurate prune technique in multi-target visiting searching problem to avoid searches on infeasible paths.

## TRAJECTORY OPTIMIZATION AND DATASET GENERATION

### Sims-Flanagan Transcription

The low-thrust (LT) trajectory between two asteroids is modelled using Sims–Flanagan Transcription (SFT)[14, 15]. As shown in Figure 1, trajectory is discretized into equal-time segments with bounded impulses applied at the segment centers as control points. The spacecraft sate $\boldsymbol{s}$ is given by:

$$\boldsymbol{s} = [\boldsymbol{p} \quad \boldsymbol{v} \quad m] \tag{1}$$

where $\boldsymbol{p}$ and $\boldsymbol{v}$ are the position and velocity of the spacecraft, $m$ is the spacecraft mass. The state is propagated forward in time from the beginning asteroid and backward from the ending asteroid. The position and velocity will be propagated using Kepler equation and the spacecraft mass will be updated using the Tsiolkovsky rocket equation. The state will be updated instantaneously after the impulses. The SFT convert the trajectory optimization into an optimal control problem, where the objective is to maximize the final mass of the spacecraft. The constraints of the problem are the magnitude of the thrusts at control points and the continuity of the spacecraft state vector at the match point, given by:

$$\boldsymbol{0} \leq \boldsymbol{T} \leq \boldsymbol{T}_{max} \tag{2}$$

$$\boldsymbol{c}_{mp} = \boldsymbol{s}_{mp}^{fwd} - \boldsymbol{s}_{mp}^{bwd} = \begin{bmatrix} p_{mp}^{fwd} - p_{mp}^{bwd} \\ v_{mp}^{fwd} - v_{mp}^{bwd} \\ m_{mp}^{fwd} - m_{mp}^{bwd} \end{bmatrix} = \begin{bmatrix} \delta p_{mp} \\ \delta v_{mp} \\ \delta m_{mp} \end{bmatrix} = \boldsymbol{0} \tag{3}$$

where $\boldsymbol{T} = [T_x \; T_y \; T_z]$ is the control vector of each segment, the superscript $fwd$ and $bwd$ indicate forward propagation and backward propagation. To further improve the robustness of the optimization process, the unit control vector $\boldsymbol{u} = [u_x \; u_y \; u_z]$ is introduced[16], where

$$T_i = T_{max} \cdot u_i \quad (4)$$

the components of $\boldsymbol{u}$ are bounded in the range of [-1,1].

The SFT converts the continuous thrust optimization problem into an optimal control problem, which can be solved by using nonlinear programming (NLP) solvers [17].

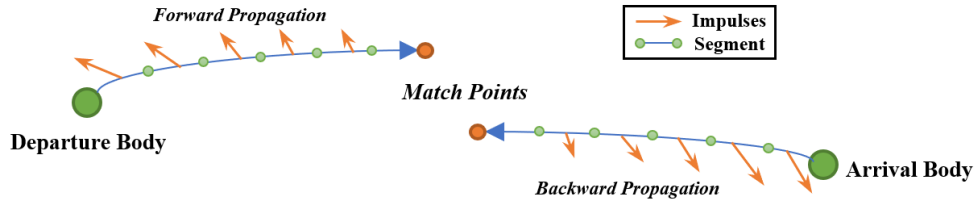

**Figure 1. The propagation in Sims-Flanagan Model. The states are propagated forward and backward to a match point, where the states must be continuous.**

Sims-Flanagan Transcription (SFT) mode has progressively been improved for higher fidelity [17, 18], and higher computation efficiency [16, 19, 20]. It is also used in mission design tools such as Mission Analysis Low-Thrust Optimization (MALTO) and Evolutionary Mission Trajectory Generator (EMTG) [21, 22], etc.

According to several presentations and publications[19, 21-24], SFT model can achieve medium to medium-high fidelity, which is sufficient for constructing a LT-feasibility identifier. In addition, SFT is also very fast because it does not require numerical integration of differential equations. And the convergence of an NLP solver solving a SFT problem is robust to poor initial guesses. These features make it ideal for automated parallel data generation on the high-performance computing cluster.

**Trajectory Data Generation**

The ephemerides are obtained from JPL HORIZON system. We developed a MATLAB - based tool to automatically acquire, pack and download SPK-formatted ephemerides via the telnet interface. The ephemerides files are read and processed by MICE (The MATAB version of SPICE toolkit) on MATAB .

The low-thrust (LT) data generation parameters are summarized in Table 1. The spacecraft is assumed to have the initial mass of 3,000 kg (including 1,000 kg dry mass), equipped with a notional electric propulsion engine with a maximum thrust of 236 mN and a specific impulse of 4,000 seconds. These specifications are chosen based on the NEXT (NASA Evolutionary Xenon Thruster) engine from Glenn Research Centre [25].

**Table 1. Studying Parameters for Trajectory Generation.**

| Parameter | Value |
|---|---|
| $m_{dry}$ | 1,000 kg |
| $m_i$ | 1,000 kg – 3,000 kg |
| $I_{sp}$ | 4,190 s |
| $T_{max}$ | 236 mN |
| SFT segment | 20 |
| $\Delta t_{impls}$ | $[60, 4 \cdot 365]$ days |
| $\Delta t_{LT}$ | $[1.2 \cdot \Delta t_{impls}, \min\{2 \cdot \Delta t_{impls}, 4 \cdot 365\}]$ days |
| Mission start date | 1 Jan 2020 – 1 Jan 2060 |

Initially, an epoch and two asteroids are randomly selected, and the ephemerides of the departure body can be determined. Then the guess of time of flight $\Delta t_{impls}$, which is the flight time of the impulsive Lambert transfer with lowest delta-V, is grid-searched and calculated by solving Lambert problem. Generally, the low-thrust transfer time is no shorter than $\Delta t_{impls}$, which gives us a range for the initial guess of feasible low-thrust transfer time $\Delta t_{LT}$. The LT transfer time is then randomly selected within $\Delta t_{LT} \in [1.2 \cdot \Delta t_{impls}, \min\{2 \cdot \Delta t_{impls}, 4 \cdot 365\}]$. Given the ephemerides, $\Delta t_{impls}$ and spacecraft configurations, the trajectory is then optimized using an NLP solver. As the LT data generation process of asteroid pairs are independent to each other, the whole process is paralleled on UNSW High Performance Computing (HPC) cluster Katana.

## DNN MODELLING FOR FEASIBLE LOW-THRUST TRANSFER IDENTIFICATION

The feasible low-thrust transfer identification problem can be seen as a binary classification problem, where the two classes are feasible and infeasible LT transfer. In this section, the DNN-based classifier as well as the other four classification models will be introduced. The feature engineering, including feature construction and preliminary feature selection technique will be first discussed. The DNN structure and optimizable hyperparameters are also determined. The Bayesian Optimization process and the evaluation metrics are also briefly introduced.

### Feature Engineering

Feature engineering is essential for training machine learning models, as the generalisation ability largely depends on the quality of it. The features considered in this work not only includes the ones required in original OCPs, but also includes newly constructed features.

The basic features are the spacecraft initial mass $m_0$, time of flight $\Delta t_{LT}$ and the ephemerides of two targets $Eph$. The ephemerides could be given in different format, such as Classical Orbital Elements ($COE$), Modified Equinoctial Elements ($MEE$), Position-Velocity ($PV$) in different coordinate systems. In addition, the differences between the features of two targets, such as $\delta_{COE}$, $\delta_{MEE}$ and $\delta_{PV}$, are also considered as the input of the networks.

The constructed features include specific angular momentum $\hbar$ ($\hbar = \boldsymbol{p} \times \boldsymbol{v}$), specific orbital energy $E$ ($E = -\mu/2a$), impulsive delta-V calculated using low-thrust transfer time and the initial flight time guess $\Delta t_{impls}$. The features are summarized in Table 2, totally 103 features are considered.

**Table 2. Features Summary.**

| Feature | Introduction |
|---|---|
| $m_i$ | Initial mass |
| tof | Time of flight |
| $COE, \delta COE$ | Classical Orbital Elements, $\delta COE = COE_1 - COE_2$ |
| $MEE, \delta MEE$ | Modified Equinoctial Elements, $\delta MEE = MEE_1 - MEE_2$ |
| $PV, \delta PV$ | Position-Velocity in Cartesian coordinates, $\delta PV = PV_1 - PV_2$ |
| $Sph, \delta Sph$ | Position-Velocity in Spherical coordinates, $\delta Sph = Sph_1 - Sph_2$ |
| $Cylind, \delta Cylind$ | Position-Velocity in Cylindrical coordinates, $\delta Cylind = Cylind_1 - Cylind_2$ |
| $R, \delta R$ | Distances from sun, $\delta R = R_1 - R_2$ |
| $E, \delta E$ | Specific orbital energy, $\delta E = E_1 - E_2$ |
| $H, \delta H$ | Specific angular momentum, $\delta H = H_1 - H_2$ |
| $tof_{ini}$ | $\Delta t_{impls}$, the time-of-flight guess provided by impulsive transfer search |
| $LamdV$ | The Lambert delta-V for the given transfer configuration |

The features $x$ will need to be standardized before training, which is given by:

$$x'_{Std} = \frac{x - \mu_x}{\sigma_x} \cdot \sigma_y + \mu_y \tag{5}$$

where $x'_{Std}$ is the standardized attributes, $x$ is the original attribute; $\mu_x$ and $\sigma_x$ are the current mean and standard deviation of the feature, and $\mu_y$ and $\sigma_y$ are the target mean and standard deviation.

The feature importance is estimated by using the tree ensembles embedded method. The embedded method earns feature importance as part of the ensemble learning process. Once the model is trained, the feature importance can be obtained.

The ensemble tree model is firstly trained on 4,000 samples. For each tree learner in the tree ensemble, the feature importance is estimated by summing changes in the risk due to splits on every predictor and dividing the sum by the number of branch nodes. The change in the node risk is the difference between the risk for the parent node and the total risk for the two children, thus for the split predictor, the increment on the importance is given by:

$$\Delta Rsk = \frac{(Rsk_p - Rsk_{child1} - Rsk_{child2})}{N_{branch}} \tag{6}$$

where $Rsk_p$ and $Rsk_{child1}$ are the node risks of parent node and children node, and $N_{branch}$ is the total number of branch nodes.

The final feature importance of ensemble tree is estimated by calculating the weighted average of feature importance score from each tree learners. The sorted feature importance rankings are given in the Figure 2, where the first sixty important features are given. The two highest ranked features are the time of flight and initial mass, both of which are needed in the original OCPs. Meanwhile, two of the constructed features, initial time-of-flight guess and Lambert delta-V are also highly relevant to the model representation. The $PV$ in different coordinates are found to be less important than $MEE$ and $COE$. The specific angular momentum is more important than Specific orbital energy.

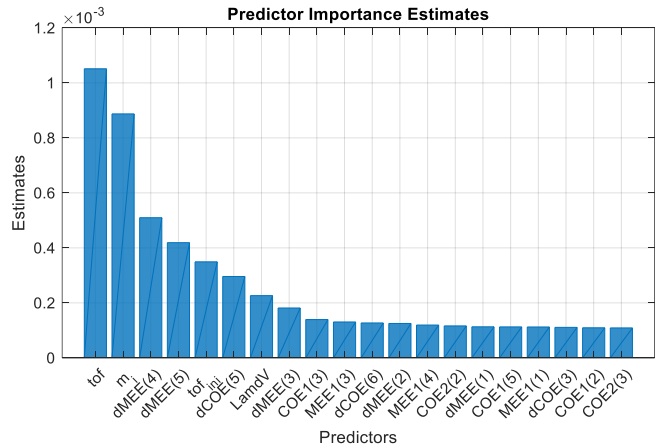

(a) feature importance of 1-20 features

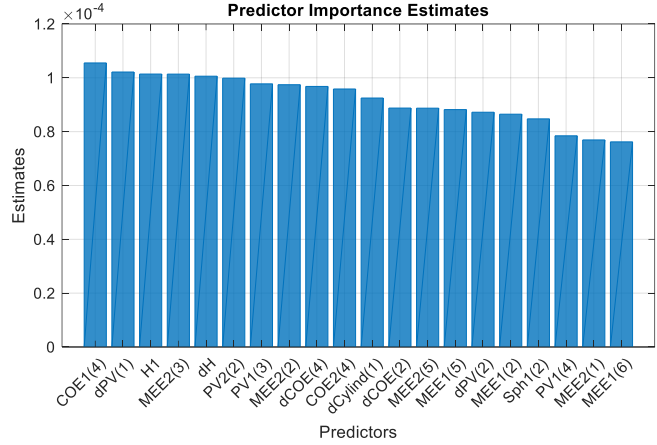

(b) feature importance of 21-40 features

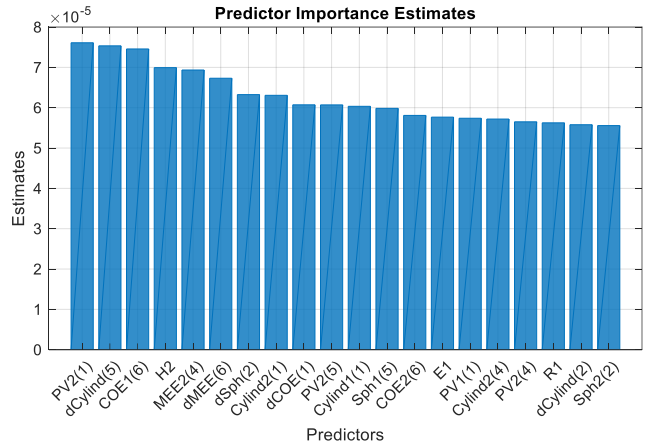

(c) feature importance of 41-60 features

**Figure 2. Sorted feature importance bar plot of 60 features. Note that plot (a), (b) and (c) have different scales. The details of the features can be found in Table 2**

**DNN-Classifier Structure**

The constructed classifier has a series deep neural network structure, where layers are arranged one after the other. The DNN contains $L$ fully connected layers, each of which has $n_{neuron}$ neurons and is followed by an activation layer.

The inputs are the standardized features introduced in the above section, the dimension of the input layer depend on the number of selected features. The outputs are probability of two classes representing the predictions of feasibility of LT transfers.

The model will be trained batch by batch. The objective of learning is to minimize the loss on mini-batches, which is defined using Cross-entropy:

$$loss_{cl} = -\sum_{i=1}^{N_b}\sum_{j=1}^{N_c} l_{ij} \ln p_{ij} \tag{7}$$

where $N_b$ is the number of samples in current batch, $N_c = 2$ is the number of classes, $l_{ij}$ is the indicator that the $i_{th}$ sample belong to class $j$, and $p_{ij}$ it is the probability that the DNN classifier associates the $i_{th}$ input with class $j$.

The optimizer parameters are updated using Adam (adaptive moment estimation), which keeps an element-wise moving average of both the parameter gradients $m_t$ and their squared values $v_t$ [26]:

$$\begin{cases} m_t = \beta_1 m_{t-1} + (1-\beta_1)\nabla_{\theta_t}\mathcal{L}(F(x_{batch}, y_{batch})) \\ v_t = \beta_2 v_{t-1} + (1-\beta_2)[\nabla_{\theta_t}\mathcal{L}(F(x_{batch}, y_{batch}))]^2 \end{cases} \tag{8}$$

where $\beta_1$ and $\beta_2$ are the gradient decay factor and squared gradient decay factor, $\mathcal{L}(F(x_{batch}, y_{batch}))$ is the loss of current batch of data, $\theta$ is the network parameter. To counteract $m_t$ and $v_t$ being biased towards zeros [26], the bias-corrected estimates are used:

$$\begin{cases} \hat{m}_t = \dfrac{m_t}{1-\beta_1^t} \\ \hat{v}_t = \dfrac{v_t}{1-\beta_2^t} \end{cases} \tag{9}$$

ADAM uses the moving averages to update the network parameter $\theta$ by:

$$\theta_{t+1}^l = \theta_t^l - \eta \frac{\hat{m}_t}{\sqrt{\hat{v}_t} + \epsilon} \tag{10}$$

where $\epsilon$ is a small constant added to avoid division by zero.

**Hyperparameter Optimization**

The hyperparameters introduced in the above section needs to be optimized, such that the model can achieve a high prediction accuracy. Hyper-parameters that need to be searched include number of fully connected layers $L$, number of neurons in each layer $n_{neuron}$, activation function $\sigma_{act}$, gradient decay factor $\beta_1$, squared gradient decay factor $\beta_2$, batch size $N_b$, initial learning rate $\eta$, learning rate drop period $T_{drop}$, learning rate drop factor $d_\eta$.

The search parameters and corresponding search spaces are defined in Table 3

**Table 3. DNN-classifier search space.**

| Parameters | Search Space |
| --- | --- |
| $L$ | [5,10] |
| $n_{neuron}$ | [200,500] |
| $\sigma_{act}$ | ReLU, leaky ReLU , ELU |
| $\beta_1$ | [0.85, 0.95] |
| $\beta_2$ | [0.9, 0.999] |
| $N_b$ | [200,800] |
| $\eta$ | [0.001,0.01] |

| $T_{drop}$ | [2, 8] epoch |
| $d_\eta$ | [0, 0.8] |

The objective of hyperparameter optimization is to find a set of hyperparameters that maximize the prediction accuracy. However, training DNN models based on given sets of hyperparameters is time consuming as a large amount of data need to be processed repeatedly. In addition, the hyper parameters presented in Table 3 result in a large search space, and various types of optimizable parameter (continuous, integer, categorical) need to be considered. To deal with these issues, Bayesian Optimization (BO) method is adapted. BO process consists of four components, including a hyperparameter search space to sample from, an objective function for hyperparameter sets evaluation, statistical model (also known as surrogate model) for modelling the objective function, and an acquisition function to maximize to determine the next configuration for evaluation[27]. Compared with the control variates method used in the literature[3, 5, 7], the BO is a sequential model-based optimization (SMBO) and can intelligently configure parameters toward a global optimum in a minimum number of steps based on the previous training history [28]. To avoid local minima, we repeatedly performed BO process with a maximum 30 trails in each process. The optimized DNN models can be found in the simulation section.

**Evaluation Metrics**

The evaluation metrics include accuracy, recall, precision, F-measure, ROC (Receiver Operating Characteristic) curve and AUC (Area Under the ROC Curve).

Accuracy is the most intuitive performance measure which is the ratio of correctly predicted observation to the total observations, that is,

$$Acc = \frac{n_{TP} + n_{TN}}{n_{TP} + n_{TN} + n_{FP} + n_{FN}} \quad (11)$$

where $n_{TP}$ and $n_{TN}$ are the numbers of true positive (TP) and true negative (TN) samples; $n_{FP}$ and $n_{FN}$ are the number of false positive (FP) and false negative (FN) samples.

The precision $Prec$ and recall $Rec$ are defined as:

$$\begin{cases} Prec = \dfrac{n_{TP}}{n_{TP} + n_{FP}} \\ Rec = \dfrac{n_{TP}}{n_{TP} + n_{FN}} \end{cases} \quad (12)$$

Recall expresses the ability to find all LT feasible samples in the dataset, precision expresses the proportion of the LT samples the classifier says were feasible actually were feasible. The F-measure is used to balance the precision and recall, and this metric is used as the replacement of prediction accuracy in the Bayesian Optimization process when dealing with imbalanced LT trajectory dataset.

The common form can be expressed as:

$$F_k = (1 + k^2) \frac{Prec \cdot Rec}{k^2 \cdot Prec + Rec} \quad (13)$$

where we consider recall is $k$ times more important than precision.

The ROC curve is also introduced for evaluating binary DNN-classifiers, the probability graph shows the True Positive Rate (TPR) against False Positive Rate (FPR) at various threshold values. The AUC is the aera under the ROC curve, which can be calculated by:

$$AUC = \sum_{i \in (TP+FP+FN+TN)} \frac{(TPR_i + TPR_{i-1}) \cdot (FPR_i - FPR_{i-1})}{2} \quad (14)$$

where $TPR = n_{TP}/(n_{TP} + n_{FN})$ and $FPR = n_{FP}/(n_{FP} + n_{TN})$. The AUC is equal to the probability that the classifier will rank a randomly chosen convergence sample higher than a randomly chosen divergence sample.

**Data Augmentation**

Based on the statistical results of our testing, when selecting two targets completely randomly, the low-thrust (LT) optimization convergence rate is around 7%. If limiting the differences between some orbital elements of two targets, such as eccentricity and inclination, the LT optimization convergence percentage is around 20%. This indicates that the raw datasets for training the classification model is heavily biased towards infeasible LT samples. In another word, using the same number of feasible LT samples and infeasible LT samples will leave a large amount infeasible samples unused. Note that, by the feasibility of a LT transfer, we mean the convergence of that LT trajectory optimization.

It will be shown in the simulation section that training on a skewed database leads to a degradation on prediction accuracy. Therefore, what we expect is to utilize more samples without influencing the balance of training set. Therefore, data augmentation technique ADAptive SYNthetic (ADASYN) [29, 30] is introduced to oversample feasible LT datasets, such that the dataset can be re-balanced, and the models can be trained with less samples.

ADASYN considers the distribution of density, and it adaptively shifting the classification decision boundary toward the difficult examples. Other oversampling might also work well in the scenario in this work, but the tuning and testing oversampling technique is beyond the scope of this work.

**SIMULATION RESULTS**

The experiment section tests the requirements of training dataset sizes for different models, perform and validate the feature selection, analyze the data augmentation quality and finally present optimized classifiers.

The performance of DNN-classifier is compared with other four representative classification models, including the Ensemble learning model (ensembled trees), K-Nearest Neighbors (KNN) model, Support Vector Machine (SVM) and Tree based model.

**Requirements on Training Dataset Sizes**

Different classification models have different requirements on training dataset sizes. This experiment provides a training set size requirement curve for the five tested models, which can be used as a guidance to choose training set sizes for the models to achieve the desirable prediction accuracies.

The performances of classifiers on difference sizes of training sets are evaluated and the results are illustrated in Figure 3. In Figure 3 (a), the accuracy of each point is obtained by averaging the results of five independent computations. In each computation, the training datasets are constructed randomly based on the given sizes, and the hyperparameters are optimized by the Bayesian Optimization (BO). For methods in Figure 3 (a), the computation became extremely slow when training on large datasets (training set size >= 30,000). For efficiency consideration, only critical parameters are optimized. In Figure 3 (b), each DNN accuracy result is the average of two computations. Note that in this experiment, datasets used are balanced.

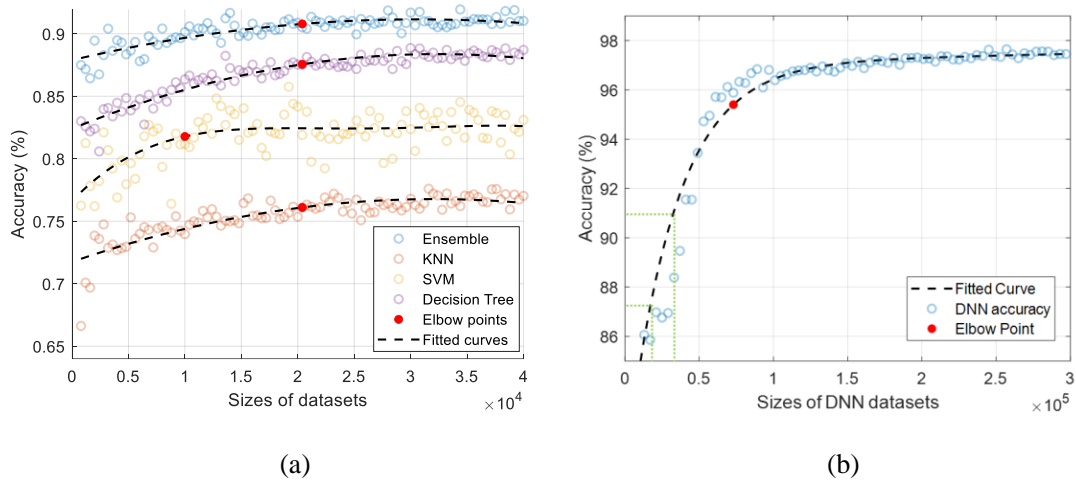

(a) (b)

**Figure 3. The performances of machine learning models and DNN model on different sizes of balanced training dataset. The elbow points indicate the points after which the datasets sizes will not increase the accuracy significantly. (a) Accuracy of Ensembled tree, KNN, SVM and tree-based method on various sizes of datasets. (b) Accuracy of the DNN model on different dataset sizes.**

It can be observed from Figure 3 (a) that, among the four machine learning methods, the ensemble model achieves the highest accuracy on all given dataset sizes. The elbow points can be seen as the minimum data amounts needed by the models to reach its generalization limits. Both the ensemble method, KNN and decision tree have an elbow point of around 20,000, while the SVM only required 10,000 samples to reach the elbow point.

By contrast, the DNN model generally required more samples to achieve expected accuracies. Given 20,000 samples, the ensemble methods can achieve and accuracy of 91%, while the DNN model provides an accuracy of 87%. To achieve the same accuracy as ensemble methods (91%), the DNN model would require over 30,000 samples. On the other hand, as it can be seen from Figure 3 (b), the DNN model is able to provide a much higher prediction accuracy when larger amount of data is available.

The data amount requirements curves presented in this section can be used as a reference for the trade-off between prediction accuracy and training cost. If pursuing a high prediction accuracy on feasible LT transfer identification, DNN-classifier should be selected. If medium accuracy is acceptable, then both the ensemble model and DNN model are workable. However, the DNN-classifier has the potential to achieve a higher accuracy – it can be trained and improved to achieve and accuracy of approximately 98% when more trajectory data are available, while the prediction accuracy of ensemble methods is limited to approximately 91%.

**Reducing Data Requirement by Using Data Augmentation Technique**

In this experiment, the influence of imbalance training dataset is evaluated first. The models are trained using databases containing 10,000 feasible low-thrust (LT) samples (as minority) and various sizes of infeasible LT samples. The chosen number of feasible LT samples is based on the data amount requirement results provided in the previous section. The accuracies, precisions and recalls of trainings on biased datasets are shown in Figure 4 (a) - Figure 6 (a). It can be observed that when adding more infeasible LT samples, the precisions of all methods generally increase, but the recalls drop dramatically. The reason why precisions become higher is because the large amount of infeasible LT samples involved in the training enables the classifiers to identify majority class more accurately, thus the False Positive number decreases. This leads to an increment on precision. On the other hand, the model has a weaker ability to identify feasible LT samples, which leads to a lower recall. It can be concluded that simply introducing more infeasible samples into the training set will cause degradation on prediction performance.

The models are then retrained using the ADASYN-augmented training data, where the feasible LT and infeasible LT samples are re-balanced. The results are shown in Figure 4 (b) - Figure 6 (b). Comparing plots (a) and (b), we see that after the augmentation, the precisions slightly decrease, and the recalls have a huge improvement.

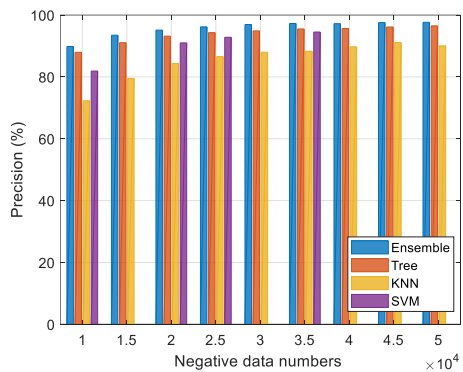
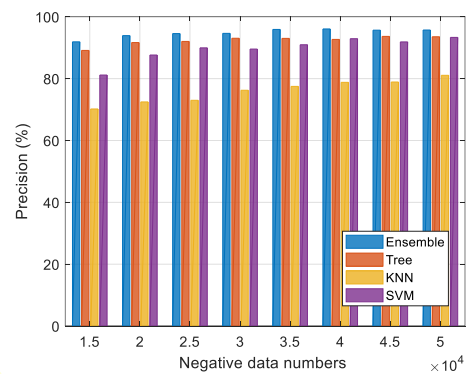

(a) Precisions on biased datasets　　　　　　　(b) Precisions on augmented datasets

**Figure 4. Precisions of machine learning methods on biased and augmented datasets.**

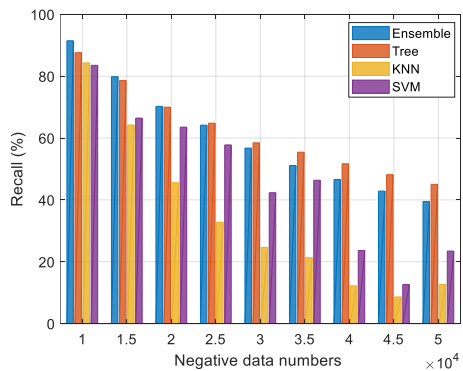
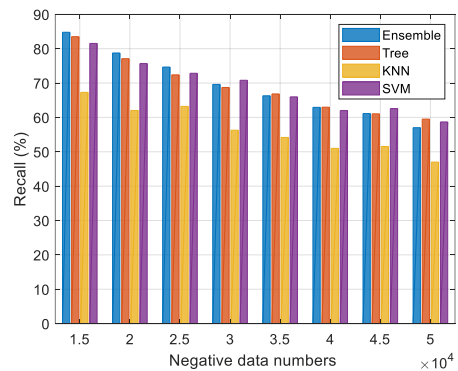

(a) Recalls on biased datasets　　　　　　　(b) Recalls on augmented datasets

**Figure 5. Recalls of machine learning methods on biased and augmented datasets.**

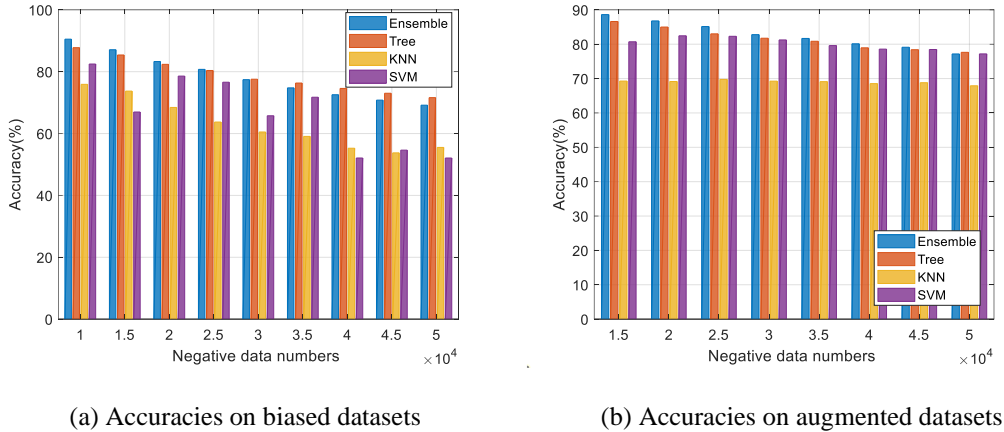

(a) Accuracies on biased datasets  (b) Accuracies on augmented datasets

**Figure 6. Accuracies of machine learning methods on biased and augmented datasets.**

For DNN-classifier evaluation, the biased training dataset is set to have 50,000 feasible LT samples and 250,000 infeasible LT samples. To reduce the level of recall degradation, the objective of model optimization is set to be $F_{10}$ (F-measure, $k = 10$), where recall is considered to be ten times more important than precision. The performances of DNN classifier on biased, augmented and balanced datasets are presented in Figure 7. As illustrated in the figure, the recall decreases by about 6% on biased training dataset, but this degradation becomes negligible after augmentation. Meanwhile, both accuracy and precision drop by approximately 1% after data augmentation.

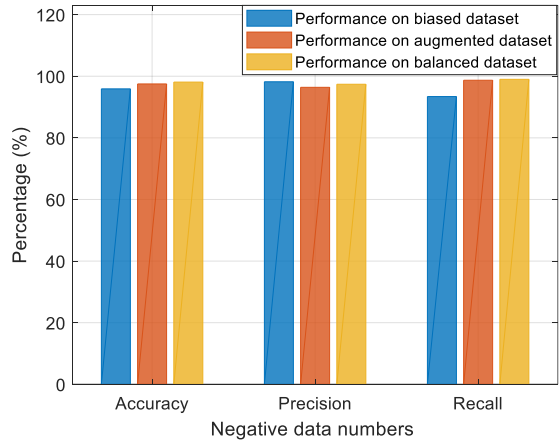

**Figure 7. Performance of DNN-classifier on balanced, imbalanced and augmented datasets.**

To conclude, although the results validate the effectiveness of data augmentation, for both the machine leaning models and DNN models, the degradation on recall and overall accuracy is inevitable as the augmented data will always introduce bias into the database. More advanced data augmentation technique and parameter tuning may reduce this bias to lower levels.

**Feature Selection**

The feature importance is analyzed in the Feature Engineering section, but the number of features needed for classification is yet to be determined. In this section, DNN classifiers are trained and evaluated by using increasing number of features (features were priorly sorted according to the feature importance analyses). For each feature set, the DNN-model is trained and then optimized by Bayesian Optimization. At this stage, to reduce the hyper-parameter search space and accelerate the optimization, only the parameters related to the DNN structures are optimized, including the number of layers, number of neurons in each layer and activation layer types. The results are shown in Table 4.

**Table 4. DNN-classifiers performance on different feature sets.**

| Feature Num. | Accuracy | Precision | Recall | F1score | F10score | AUC |
|---|---|---|---|---|---|---|
| 30 | 0.9622 | 0.9495 | 0.9763 | 0.9627 | 0.9760 | 0.9922 |
| 40 | 0.9739 | 0.9661 | 0.9824 | 0.9741 | 0.9822 | 0.9950 |
| 50 | 0.9735 | 0.9651 | 0.9825 | 0.9737 | 0.9823 | 0.9951 |
| 55 | 0.9744 | 0.9659 | 0.9835 | 0.9746 | 0.9833 | 0.9951 |
| 60 | 0.9751 | 0.9668 | 0.9840 | 0.9753 | 0.9839 | 0.9953 |
| 65 | 0.9741 | 0.9668 | 0.9818 | 0.9743 | 0.9817 | 0.9950 |
| 70 | 0.9730 | 0.9654 | 0.9813 | 0.9732 | 0.9811 | 0.9948 |
| 103 | 0.9522 | 0.9430 | 0.9625 | 0.9526 | 0.9623 | 0.9864 |

The performances of DNN models are close when given first 40 - 70 features. Among these seven models, we can see that the models trained by using first 60 features have the best performances. Thus, the first 60 features are selected as new feature set for training machine learning and DNN models. The detailed selected feature elements can be found in Figure 2.

To validate the feature selection, the training processes of DNN models trained by using all features and selected features are compared in the Figure 8. As shown in the plot, the training process after feature selection is more stable and the validation curve eventually reaches a higher accuracy. The detailed evaluation metrics comparison of two DNN models can be found in the next section.

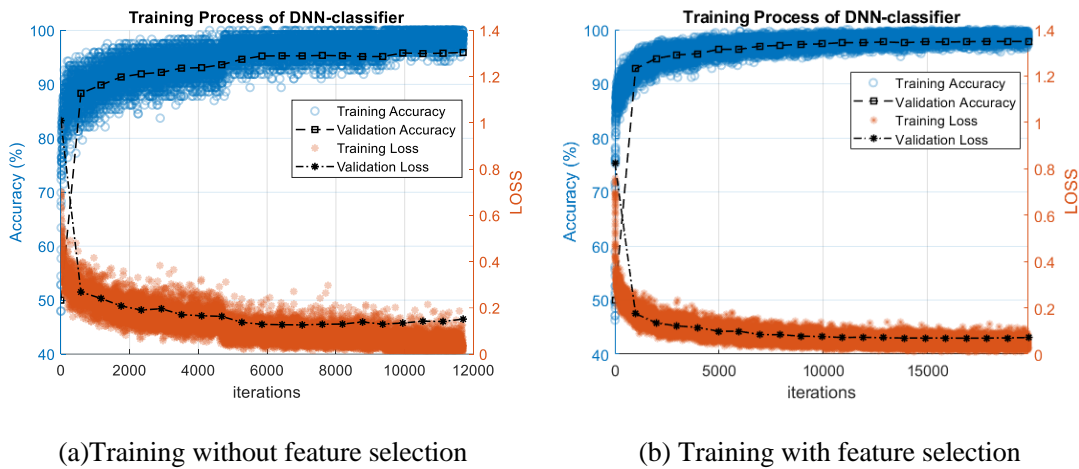

(a) Training without feature selection      (b) Training with feature selection

**Figure 8. Training progresses of two DNN-models. Each point represents the evaluation result on a mini batch. The models are validated every epoch.**

It should be noted that, although the feature selection is proofed to be effective, the selected features may still be redundant. More work on DNN model interpretation should be performed in the feature work.

**Final Evaluation of the DNN-Classifier and Machine Learning Models**

In the above sections, we have determined the data amount requirements and feature configurations for low-thrust (LT) feasibility identification. These knowledges will be used for training machine learning models and DNN-classifiers in this section.

The training of machine learning models uses 50,000 samples (except for SVM, we used only 30,000 samples for efficiency consideration). The DNN training uses 400,000 samples. The hyperparameter optimization is performed by using parallel Bayesian Optimization (BO), where all relevant hyperparameters are considered and optimized. The objective function evaluation process is provided in the Figure 9 (a) - (e). As can be seen from these optimization process, the optimal hyperparameter configurations are found within 30 evaluations. To accelerate the training and optimization process, the BO is performed parallelly on UNSW High Performance Computing (HPC) cluster Katana. As shown in Figure 9 (f), the objective evaluation time for DNN-classifier training is reduced from 17 hour (sequentially) to 2.7 hours (parallelly).

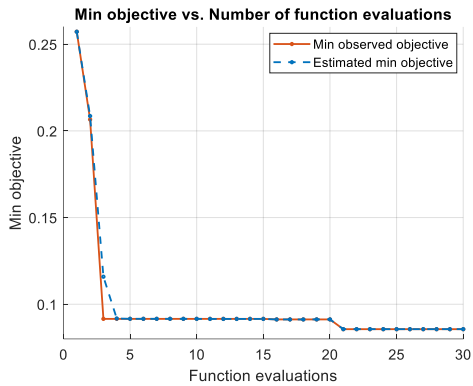

(a) Ensemble Model

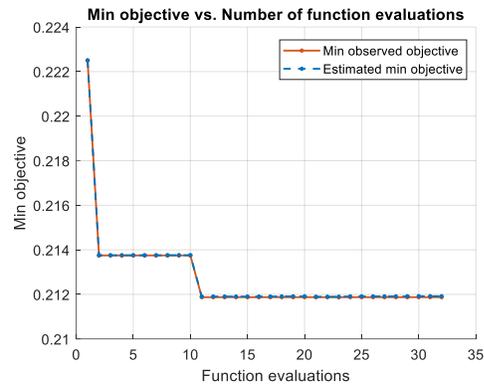

(b) KNN model

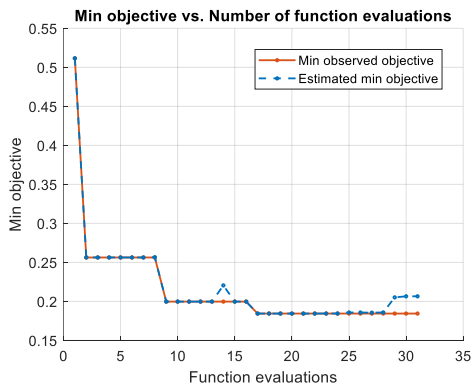

(c) SVM

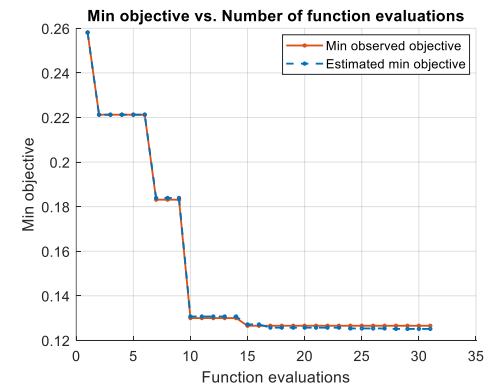

(d) Tree based model

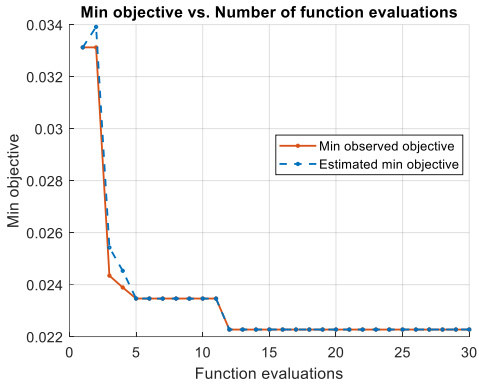
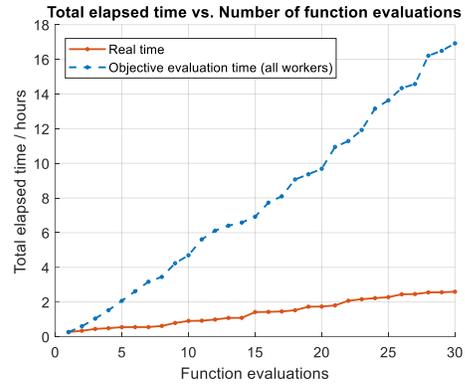

(e) DNN-classifier  (f) DNN-classifier evaluation time

Figure 9. Bayesian optimization process of classifiers.

The optimized hyperparameters are summarized in Table 6. Two optimized hyperparameter configurations both use Leaky ReLU as activation layer and the initial learning rate is very close. However, the DNN-classifier trained using selected features allows deeper and wider network structure.

Table 5. The optimized hyperparameter for two DNN classifiers.

| Parameters | DNN1(with feature selection) | DNN2(without feature selection) |
|---|---|---|
| $L$ | 9 | 6 |
| $n_{neuron}$ | 440 | 266 |
| $\sigma_{act}$ | Leaky ReLU | Leaky ReLU |
| $\beta_1$ | 0.9041 | 0.9230 |
| $\beta_2$ | 0.9540 | 0.9711 |
| $N_b$ | 500 | 123 |
| $\eta$ | 0.0031 | 0.0031 |
| $T_{drop}$ | 6 | 8 |
| $d_\eta$ | 0.1720 | 0.3284 |

The performances of classifiers are further given in Table 6 and Figure 10. The DNN-classifier trained using selected features showed the best performance on all evaluation metrics, and ensemble models shows the best performance among conventional classifiers. By comparing two DNN models, we can conclude that the feature selection technique is effective in improving the DNN model performance. In this case, feature selection improves the DNN prediction accuracy, precision and recall by 2.7%, 2.8% and 2.5%, respectively. As illustrated in Figure 10, the DNN classifiers curve is above other models, which also indicate the DNN classifier has a better representation.

Table 6. The performances of optimized classifiers.

| Model | Accuracy | Precision | Recall | F1score | F10score | AUC |
|---|---|---|---|---|---|---|
| Ensemble | 0.9108 | 0.8967 | 0.9285 | 0.9123 | 0.9282 | 0.9649 |
| KNN | 0.8018 | 0.7565 | 0.8900 | 0.8178 | 0.8884 | 0.8823 |
| SVM | 0.8385 | 0.8348 | 0.8440 | 0.8394 | 0.8439 | 0.9135 |
| Tree | 0.8723 | 0.8751 | 0.8685 | 0.8718 | 0.8686 | 0.9427 |
| DNN1 | 0.9790 | 0.9707 | 0.9877 | 0.9792 | 0.9876 | 0.9953 |
| DNN2 | 0.9522 | 0.9430 | 0.9625 | 0.9526 | 0.9623 | 0.9864 |

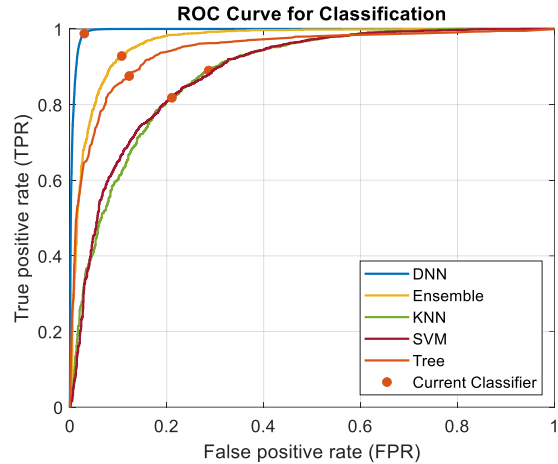

**Figure 10. ROC curves of multiple classifiers from Table 6. The closer the ROC curve get to the upper left corner, the better the test is overall. The red dot is the performance of current classifier.**

In this experiment, we used balanced dataset for training, and the optimization objective is set to maximize the overall prediction accuracy. In fact, the optimization objective can also be set as F-measure ($k > 1$) in the feasible LT trajectory generation scenario, because we prefer high recall over high accuracy/precision. Recall represents the ability to identify feasible LT samples, and a low recall will lead to a biased feasible LT database as many feasible LT samples are classified as infeasible. By contrast, lower precision brings more false-positive (FP) samples, which only leads to a slightly higher database generation time. However, for other application scenarios, such as multi-targets visiting sequence search, where precision and recall are equally important, the accuracy or F1 should be used as the optimization objective.

## CONCLUSION

In this study, we present the construction, training and optimization a DNN-classifier that can identify feasible low-thrust transfer with an accuracy of 97.9% and a recall of 98.8%. Comparing DNN-classifiers with other four representative machine learning models, we found that the DNN-classifier is the most advanced model for high-accuracy LT feasibility identification.

## ACKNOWLEDGMENTS

Financial support for this research was provided by the UIPA scholarship from University of New South Wales and CSIRO Top-up Scholarship. This research includes computations using the computational cluster Katana supported by Research Technology Services at UNSW Sydney.